\documentclass{article}
\usepackage{geometry}
\usepackage{amsmath} 
\usepackage{amssymb} 
\usepackage{MnSymbol}
\usepackage{extarrows}
\usepackage{enumerate} 

\usepackage{bm} 

\usepackage{graphicx} 
\usepackage{subfigure} 
\usepackage{multirow} 

\usepackage{cite}
\usepackage[colorlinks = true, citecolor = blue,  urlcolor = red]{hyperref}
\newtheorem{theorem}{Theorem}[section]
\newtheorem{definition}{Definition}[section]

\newtheorem{lemma}{Lemma}[section]

\newtheorem{proposition}{Proposition}[section]
\newtheorem{corollary}{Corollary}[section]
\newtheorem{remark}{Remark}[section]

\newcommand{\E}{\mathbb{E}}
\newcommand{\V}{\text{Var}}
\newcommand{\pf}{\textbf{Proof: }}
\newcommand{\e}{\hfill$\blacksquare$}

\newcommand{\R}{\mathbb{R}}
\newcommand{\KL}{\text{KL}}

\allowdisplaybreaks[4] 
\numberwithin{equation}{section} 
\title{Entropy Jump and Entropic Central Limit Theorem for Independent Sum}
\date{}
\author{Liuquan Yao~~$\cdot$~~Shuai Yuan}
\begin{document}

	\maketitle

	\noindent\textbf{Abstruct} 
	It is a manuscript for   results about entropic central limit theorem for independent sum under finite Poincar$\acute{\mbox{e}}$ constant conditions.



\section{Introduction}

\begin{definition}
	Consider a random variable $X$ with density $f$, the differential entropy of $X$ is
	$$h(X)=h(f):=-\int_\R f(x)\ln f(x) dx.$$
\end{definition}
\begin{definition}
	Set $X$ is a random variable with absolutely continuous density $f$ and its Radon-Nikodym derivative $f'$, call
	$$J(X)=J(f):=\int_{ \{f(x)>0\}} \frac{f'(x)^2}{f(x)}dx$$
	the Fisher information of $X$.
	
	If a density $\hat{f}$ satisfies $\hat{f}(x)=f(x), a.s.$ then we put
	$$J(\hat{f})=J(f).$$
	In other cases, we set $J(X)=\infty.$
\end{definition}

There is a inequility focused by some studies,  the Entropy Power Inequility(EPI)
$$e^{2h(X+Y)}\ge e^{2h(X)}+e^{2h(Y)},$$
or its corollary, the Entropy Jump Inequility(EJI)
\begin{equation}\label{EJI}
h(\frac{X+Y}{\sqrt{2}})\ge \dfrac{h(X)+h(Y)}{2},
\end{equation}
for any independent random variables $X, Y$ with densities. Carlen(1991)\cite{Carlen1991} used the  Ornstein-Uhlenbeck Semigroup to prove the existance of the non-zero gap $\delta$ in EJI and the entropic central limit theorem can be a directly corollary. 

Johnson(2000)\cite{EI-CLT} extended Carlen's work to the independent case with some special conditions(not just i.i.d.), and Ball etc.(2003)\cite{gap-EJ}, Johnson etc.(2004)\cite{FI-CLT} gave a feasible expression of the gap $\delta$ for i.i.d. case respectively by different methods. Particularly, the paper \cite{FI-CLT}  showed the speed of convergence in standard Fisher information sense is $O(\frac{1}{n}).$ Ball ect.(2004)\cite{Ball2004} also gave the speed and considered the weighted sum case. We analogize the porjection method in Johnson(2004) \cite{FI-CLT} to the indepenednt case(not i.i.d.), and obtain that for two nonsingular random variables $X, Y$ with absolutely continuous densities and their Poincar$\acute{\mbox{e}}$ constant satisfies $R^*_X, R^*_Y<\infty$, then $\exists a=a(\sigma^2_x, \sigma^2_Y), R=R(R_X^*, R_Y^*)>0$,
	$$h(\dfrac{X+Y}{\sqrt{2}})-\dfrac{h(X)+h(Y)}{2}\ge \dfrac{\min\{a,1\}}{\min\{a,1\}+2R}\left[  \dfrac{1}{2}\ln(2\pi e\dfrac{\sigma_X^2+\sigma_Y^2}{2})- \dfrac{h(X)+h(Y)}{2}\right].$$

Recently,  Johnson(2020)\cite{Johnson2020} relaxed the finite Poincar$\acute{\mbox{e}}$ constant condition for entropy gap theory used some functional analysis tools.


\section{Fisher Information Inequalities for Convolution}
The convolution Fisher information inequality is essential theory for entropy jump and entropic CLT. Paper \cite{FI-CLT} gave a convolution Fisher information inequality for i.i.d. case under finite Poincar$\acute{\mbox{e}}$ constant condition
.
	\begin{definition}\label{Poincare}
	Give a random variable $X$, define
	$$R_X^*=\sup_{g\in H_1^*(X)}\dfrac{\E[g(X)^2]}{\E[g'(X)^2]},$$
	where $H_1^*(X):=\{ f| f~is~absolutely~continuous, \V[f(X)]>0, \E[f(X)]=0, \E[f(X)^2]<\infty, \E[f'(X)]=0  \}.$ $R^*_X$ is called the Poincar$\acute{\mbox{e}}$ constant of $X$.
\end{definition}

By the definition, $R^*_{aX}=a^2R^*_X,\;\forall a\neq 0.$ Suppose finite Poincar$\acute{\mbox{e}}$ constant, then Johnson(2004) \cite{FI-CLT} obtained:

\begin{proposition}( \cite{FI-CLT} )\label{convolution of J in iid}
	
	Consider $X_1, X_2$ IID with absolutely continuous densities, variance $\sigma^2$ and restricted Poincar$\acute{\mbox{e}}$ constant $R^*$, then
	$$J(\dfrac{X_1+X_2}{\sqrt{2}})\le \dfrac{2R^*}{\sigma^2+2R^*}J(X_1).$$
\end{proposition}

This Proposition can be expanded to independent situation.

\begin{proposition}\label{convolution of J in independent}
	Consider two independent random variables $X, Y$, with absolutely continuous densities and finite variances $\sigma^2_X, \sigma^2_Y$, and there exists $0<R<\infty$ s.t.
	$$ R_X^*\le R, \;\; R_Y^*\le R.$$
	Denote $\dfrac{1}{2}[J(X)+J(Y)]=J_1, J(\frac{X+Y}{\sqrt{2}})=J_2, \frac{\sigma_X^2+\sigma_Y^2}{2}=\sigma^2$, then
		\begin{equation}\label{simple result like corollary}
			\dfrac{\sigma^2J_1-1}{\sigma^2+2R}\le J_1-J_2\le \dfrac{\sigma^2J_1-1}{\sigma^2}.
		\end{equation}

\end{proposition}
\pf In Section \ref{proof of independent convolution J}. \e

\section{Entropy Jump}
Using entropy to prove CLT can be trace back to 1959 \cite{Lin1959} and  entropy jump is the basic tool. Paper \cite{Carlen1991} found that entropy of convolution always increase with fixed variance unless there are all Gaussian random variables.

\begin{theorem}[\cite{Carlen1991} ]\label{CT-EJ for same variance}
	Let $X_1$ and $X_2$ be independent random variables with zero mean and $\sigma^2(X_1)=\sigma^2(X_2)=1$. Suppose
	$$J(X_1), J(X_2)\le J_0<\infty.$$
	Let $L$ be a given function on $[0, \infty)$ which decreases to zero at infinity s.t. $L_{X_1}, L_{X_2}\le L$. Also let some number $a$ satisfying $0<a<\frac{1}{2}$ be given. Suppose for some $\varepsilon>0$,
$$h(X_1)\le h(g)-\varepsilon,$$
where  $g$ is the density of standard normal distribution.
Then there is a $\delta_{\varepsilon}>0$ depending only on $J_0, L, a$ and $\varepsilon$ so that for any $\lambda$ with
$$a\le \lambda^2\le 1-a,$$
we have
$$h(\lambda X_1+(1-\lambda^2)^{1/2}X_2)-\lambda^2 h(X_1)-(1-\lambda^2)h(X_2)\ge \delta_\varepsilon.$$
\end{theorem}

For finite Poincar$\acute{\mbox{e}}$ conatant condition, such as log-concave density, specific gap for entropy jump can be deduced as a corollary of convolution Fisher information inequality.

\begin{theorem}[\cite{gap-EJ} ]
	Let $X$ be a random variable with variance 1 and finite entropy, whose density $f$ satisfies the Poincar$\acute{\mbox{e}}$  inequality
	$$R\int_\R fs^2\le \int_\R f(s')^2,$$
	for any function $s$ satisfying
	$$\int_\R fs=0,$$
	where $R>0$ is a constant.
	
	Set $Y$ is an independent copy of $X$, then
	$$h(\dfrac{X+Y}{\sqrt{2}})-h(X)\ge \dfrac{R}{2+2R}[h(g)-h(X)],$$
	where $g$ is  the density of standard normal distribution.
\end{theorem}

Finally, we give the definition of the Lebesgue adjoint to Ornstein-Uhlenbeck Semigroup which was used in our entropy jump theorem.
\begin{definition}(Lebesgue adjoint to Ornstein-Uhlenbeck Semigroup)\label{Ornstein-Uhlenbeck}
	Set $t\geq 0$,  define a map $p_t^*: \mathcal{T} \rightarrow \mathcal{T}$
	$$(p_t^*f)(x) = \int_\R f(y)g_{1-e^{-2t}}(x-e^{-t}y)dy.$$
\end{definition}
 And an important theorem, de Bruijn's identity, which shows the relationship between FIsher information and entropy.
\begin{theorem}(de Bruijn's identity\cite{cover2006elements})\label{Bruijn}
	Let $X$ be a random variable  with density  $f\in\mathcal{T}$,  $G$ is standard Gaussian distribution which is independent with $X$.  Then the function $t\mapsto h(X+\sqrt{t}G)$is continuous differentiable on $(0, \infty)$,  and
	$$
	\frac{d}{dt}h(X+\sqrt{t}G) = \frac{1}{2}J(X+\sqrt{t}G).
	$$
\end{theorem}

By de Bruijn's identity, we could deduce the following theorem from Proposition \ref{convolution of J in independent} for the bounds of entropy jump in indepdent case.
\begin{theorem}\label{entropy jump in indpendent}
		Consider two independent nonsingular random variables $X, Y$, with absolutely continuous densities and finite variances $\sigma^2_X, \sigma^2_Y$, and there exists $0<R<\infty$ s.t.
	$$  \max \{R_X^*, R_Y^*, R^*_G\}\le R.$$
	where $G$ is standard Gaussian distribution. Then there exists $c=c(\sigma^2_X, \sigma^2_Y, R)$ s.t. 
	$$c\left(\dfrac{1}{2}\ln(2\pi e\dfrac{\sigma_X^2+\sigma_Y^2}{2})- \dfrac{h(X)+h(Y)}{2}\right)\le h(\dfrac{X+Y}{\sqrt{2}})-\dfrac{h(X)+h(Y)}{2}.$$
\end{theorem}

\pf 
By de Bruijn's identity(Theorem \ref{Bruijn}),
$$h(X)=h(G)-\int_0^\infty [J(p_s^* X)-1]ds,$$
where $G$ is standard Gaussian distribution and the define of $p_s^*$ can be found in Definition \ref{Ornstein-Uhlenbeck}. So 
$$h(\dfrac{X+Y}{\sqrt{2}})-\dfrac{h(X)+h(Y)}{2}=\int_0^\infty \dfrac{1}{2}[J(p_s^* X)+J(p_s^* Y)]-J(\dfrac{p_s^* X+p_s^* Y}{\sqrt{2}}) ds.$$
Suppose $R$ is  bigger than the Poincar$\acute{\mbox{e}}$ constant of $X, Y, G$, then by the fact in \cite{gap-EJ} that:  
for any  random variables $X, Y$ satisfying $R^*_X, R^*_Y\le R$, then $\forall \lambda\in [0,1]$,
\begin{equation}\label{R not change under convolution}
R^*_{\lambda X+(1-\lambda^2)^{1/2} Y}\le R.
\end{equation}
Thus by \eqref{simple result like corollary},
\begin{align*}
&h(\dfrac{X+Y}{\sqrt{2}})-\dfrac{h(X)+h(Y)}{2}\\
&\ge \int_0^\infty \dfrac{1-(1-\sigma^2)e^{-2s}}{1-(1-\sigma^2)e^{-2s}+2R}[\dfrac{1}{2}J(p_s^* X)+\dfrac{1}{2}J(p_s^* Y)-\dfrac{1}{1-(1-\sigma^2)e^{-2s}}]ds\\
&\ge \dfrac{\min\{a,1\}}{\min\{a,1\}+2R}\left( \int_0^\infty [\dfrac{1}{2}J(p_s^* X)+\dfrac{1}{2}J(p_s^* Y) -1 ]ds  +  \int_0^\infty    \dfrac{(\sigma^2-1)e^{-2s}}{1-(1-\sigma^2)e^{-2s}} ds\right )\\
&=\dfrac{\min\{a,1\}}{\min\{a,1\}+2R}\left( \dfrac{1}{2}\ln(2\pi e)- \dfrac{h(X)+h(Y)}{2} +  \int_0^\infty    \dfrac{(\sigma^2-1)e^{-2s}}{1-(1-\sigma^2)e^{-2s}} ds\right )\\
&=\dfrac{\min\{a,1\}}{\min\{a,1\}+2R}[\dfrac{1}{2}\ln(2\pi e\sigma^2)-\dfrac{h(X)+h(Y)}{2}],
\end{align*}
where $\sigma^2=\frac{\sigma^2_X+\sigma^2_Y}{2}, a=\min \{\sigma^2_X, \sigma^2_Y  \},$ then we complete the proof by taking
$$c=\dfrac{\min\{a,1\}}{\min\{a,1\}+2R}.$$
\e
 \begin{remark}
 	If $\sigma_X, \sigma_Y\ge \sigma_0$, then we can take
 	\begin{equation}
 	c=\dfrac{\min\{\sigma^2_0,1\}}{\min\{\sigma^2_0,1\}+2R}.
 	\end{equation}
 \end{remark}

\section{Central Limit Theorem for Independent Sum}
\begin{theorem}\label{CLT for independent}
	Consider an independent nonsingular random variable sequence $X_1,X_2,\cdots$ with absolutely continuous densities, and denote the sum under the power of 2 as
	$$\hat{S}_n=\dfrac{X_1+X_2+\cdots +X_{2^n}}{2^{\frac{n}{2}}}.$$
			If the sequence $\{X_n,n\ge1\}$ satisfies the following conditions:
			\begin{enumerate}
				\item Uniformly finite Poincar$\acute{\mbox{e}}$ condition $R:=\sup_nR^*_{X_i}<\infty,\;\;\forall i\in\mathbb{Z}^+$;

				\item Stable entropy condition $$h(\dfrac{X_1+X_2\cdots+X_{2^m}}{2^{m/2}})=h(\dfrac{X_{2^m+1}+X_{2^m+2}\cdots+X_{2^{m+1}}}{2^{m/2}}),\;\;\forall m\in\mathbb{N};$$
				
				\item Convergent variance condition $\V(\hat{S}_n)\to \sigma^2>0,\; n\to\infty$, for a real number $\sigma$;
				
								\item Zero mean condition $\E(X_i)=0,\;\;\forall i\in\mathbb{Z}^+$.
				
			\end{enumerate}
		Then 
		\begin{equation}
		\|p_{\hat{S}_n}-\phi_{\sigma^2}\|_{L^1}\to 0,\;\;n\to\infty,
		\end{equation}
		where $p_{\hat{S}_n}$ is the density of $\hat{S}_n$, $\phi_{\sigma^2}$ is the density of  $\mathcal{N}(0,\sigma^2)$.
\end{theorem}
\pf  For any given $m$, by condtion 1 and \eqref{R not change under convolution}. we can use induction to show that
$$R^*_{\frac{X_1+X_2\cdots+X_{2^m}}{2^{m/2}}}\le R, R^*_{\frac{X_{2^m+1}+X_{2^m+2}\cdots+X_{2^{m+1}}}{2^{m/2}}}\le R,\;\;\forall  m\in\mathbb{N}.$$
We denote $h_m:=h(\hat{S}_m)$, then due to Theorem \ref{entropy jump in indpendent} and condition 2, we obtain
\begin{equation}\label{core inequility}
h_{m+1}-h_m\ge \dfrac{\min\{b,1\}}{\min\{b,1\}+2R}\left(  \frac{1}{2}\ln (2\pi e \V(\hat{S}_{m+1}))-h_m\right),\;\;\forall  m\in\mathbb{N},
\end{equation}
where $b:=\inf_n \V(\hat{S}_n)>0$ since the  condition 3 and the nonsingular condition hold. 

On the other hand, according to EJI \eqref{EJI}, sequence $\{h_n, n\ge 1\}$ is increasing with upper bound $\frac{1}{2}\ln(2\pi e [\sup_n \V(\hat{S}_n)])<\infty$ since the  condition 3, thus the  sequence $\{h_n, n\ge 1\}$ is convergent. This implies
$$h_{m+1}-h_m\to 0, \;\; m\to \infty.$$
Recall that $\frac{1}{2}\ln (2\pi e \V(\hat{S}_{m+1}))-h_m\ge h_{m+1}-h_m\ge 0$, from \eqref{core inequility} we obatin 
$$\frac{1}{2}\ln (2\pi e \V(\hat{S}_{m+1}))-h_m\to 0, \;\; m\to\infty,$$
therefore by condition 3,
$$\frac{1}{2}\ln (2\pi e \V(\hat{S}_{m}))-h_m\to 0, \;\; m\to\infty,$$
i.e. 
$$\KL(\hat{S}_m\| \mathcal{N}(0,\V(\hat{S}_m)) )\to 0, \;\; m\to\infty,$$
since condition 4 holds. By Pinsker inequility,
$$\|p_{\hat{S}_n}-\phi_{\V(\hat{S}_m)}\|_{L^1}\to 0,\;\;n\to\infty.$$
Then   the proof  is completed since
$$\|p_{\hat{S}_n}-\phi_{\sigma^2}\|_{L^1}\le \|p_{\hat{S}_n}-\phi_{\V(\hat{S}_n)}\|_{L^1}+ \|\phi_{\sigma^2}-\phi_{\V(\hat{S}_n)}\|_{L^1},\;\;\forall n\in \mathbb{Z}^+,$$
and condition 3.
\e

\begin{corollary}\label{speed of CLE independent}
	If the condition 3 in Theorem \ref{CLT for independent} is replaced by
	$$\V(X_i)=\sigma^2,\;\;\forall i\in\mathbb{Z}^+,$$
	then 
			\begin{equation}
	\KL(\hat{S}_n\|\mathcal{N}(0,\sigma^2))=O\left(  (\dfrac{2R}{\min\{\sigma^2,1\}+2R})^{n} \right),\;\;n\to\infty.
	\end{equation}
\end{corollary}
\pf
By \eqref{core inequility} and equal variance condition,
$$h-h_{m+1}\le(1-c)(h-h_m),\;\;\forall m\in\mathbb{N}, $$
where $h=\frac{1}{2}\ln(2\pi e\sigma^2), c=\dfrac{\min\{\sigma^2,1\}}{\min\{\sigma^2,1\}+2R}.$ Thus
$$\KL(\hat{S}_n\|\mathcal{N}(0,\sigma^2))=h-h_n\le (1-c)^n(h-h_0),$$
where $h_0=h(X_1)$. Using Pinsker inequility again we complete the proof.
\e

For normal sum $S_n:=1/\sqrt{n} \sum_{i=1}^n X_i$, we have the following two corollaries from Theorem \ref{CLT for independent}.

\begin{corollary}\label{iff theorem of Sn}
		Consider an independent nonsingular random variable sequence $X_1,X_2,\cdots$ with absolutely continuous densities, and denote the sum 
	$$S_n:=\dfrac{X_1+X_2+\cdots +X_{n}}{\sqrt{n}}.$$
	If the sequence $\{X_n,n\ge1\}$ satisfies the following conditions:
	\begin{enumerate}
		\item Uniformly finite Poincar$\acute{\mbox{e}}$ condition $R:=\sup_nR^*_{X_i}<\infty,\;\;\forall i\in\mathbb{Z}^+$;

		\item Increasing adding entropy condition $$h(\dfrac{X_1+X_2\cdots+X_{2^m}}{2^{m/2}})\le h(\dfrac{X_{2^m+1}+X_{2^m+2}\cdots+X_{2^{m+1}}}{2^{m/2}}),\;\;\forall m\in\mathbb{N};$$

		\item  Mean and variance condition $\E(X_i)=0, \V(X_i)=\sigma^2,\;\;\forall i\in\mathbb{Z}^+$.
		
	\end{enumerate}
	Then the entropy sequence $\{Ent_n:=h(S_n), n\ge 1\}$ is convergent if and only if
	\begin{equation*}
	\KL(S_n\|\mathcal{N}(0,\sigma^2))\to 0,\;\;n\to\infty.
	\end{equation*}
\end{corollary}
\pf Since $\KL(S_n\|\mathcal{N}(0,\sigma^2))=\frac{1}{2}\ln(2\pi e\sigma^2)-Ent_n$, the one part of proof is finished.

Next we suppose $\{Ent_n, n\ge1\}$ is convergent, thus its subsequence $\{h_n,n\ge 1\}$ defined in Theorem \ref{CLT for independent} is also convergent. then consider entropy jump
$h_{n+1}-\dfrac{h_n+\bar{h}_n}{2},$
where
$$\bar{h}_n=h(\dfrac{X_{2^n+1}+X_{2^n+2}\cdots+X_{2^{n+1}}}{2^{n/2}}).$$
For one hand, by Theorem \ref{entropy jump in indpendent},
$$h_{n+1}-\dfrac{h_n+\bar{h}_n}{2}\ge 0,$$
thus
\begin{equation}\label{liminf}
\liminf_{n\to\infty} h_{n+1}-\dfrac{h_n+\bar{h}_n}{2}=\liminf_{n\to\infty} \dfrac{h_{n}-\bar{h}_n}{2}\ge 0.
\end{equation}
One the other hand, by condition 2,
\begin{equation}\label{limsup}
h_{n}-\bar{h}_n\le 0, \forall n\in\mathbb{N}^+. 
\end{equation} 
Combining \eqref{liminf} and \eqref{limsup}, we have
$$\lim_{n\to\infty} h_{n}-\bar{h}_n=0,$$
which implies
$$\lim_{n\to\infty} h_{n+1}-\dfrac{h_n+\bar{h}_n}{2}=0.$$
Use Theorem \ref{entropy jump in indpendent} again we obtain
$$\lim_{n\to\infty} \dfrac{1}{2}\ln(2\pi e\sigma^2)-h_n=\lim_{n\to\infty} \dfrac{1}{2}\ln(2\pi e\sigma^2)-\dfrac{h_n+\bar{h}_n}{2}=0,$$
i.e. a subsequence of $\{Ent_n, n\ge1\}$ converge to $\dfrac{1}{2}\ln(2\pi e\sigma^2)$, then we complete the proof due to the convergence of the whole sequence  $\{Ent_n, n\ge1\}$.

\begin{corollary}\label{rate theorem of Sn for CLT}
Under the conditions in Corollary \ref{speed of CLE independent}, we have
\begin{equation}
\KL(S_n\|\mathcal{N}(0,\sigma^2))=\begin{cases}O(n^{\log c}), & \mbox{if~} c>1/2,\\
O(n/\log n), & \mbox{if~} c=1/2,\\
O(n^{-1}), & \mbox{if~} c<1/2,
\end{cases}
\;\;\mbox{as~} n\to \infty,
\end{equation}

where
$$c=\dfrac{2R}{\min\{\sigma^2,1\}+2R}.$$

\pf Denote $D_n=\KL(S_n\|\mathcal{N}(0,\sigma^2))$. Since $D_n$ is positive and subadditive, we have
\begin{equation}\label{Dn subadd}
nD_n\le 2^kD_{2^k}+2^{k-1}D_{2^{k-1}}+\cdots +D_1,
\end{equation}
where $k=\lfloor \log n \rfloor.$ By the definition of $c$ and Corollary \ref{speed of CLE independent}, when $c\neq 1/2$ we have
$$nD_n\le |\frac{1-2^kc^k}{1-2c}|.$$
Use $2^k\le n$, we obatin
$$D_n=O(n^\alpha),$$
where
$$\alpha=\max\{-1, \log c\}.$$
When $c=1/2$, then by \eqref{Dn subadd},
$$nD_n\le k+1\le \log n+1.$$
\end{corollary}

\section{ Proof of Proposition \ref{convolution of J in independent} }\label{proof of independent convolution J}

We would prove the result for Fisher information jump in independent case, following the method in \cite{gap-EJ}.

\subsection{Preparation}
This part we show some conclusions given by \cite{gap-EJ}, which still hold in just independent case. We miss the proofs since they are the same as the i.i.d. case.

	The first  is about the  projection lemma, which is the core of the proof. 
\begin{lemma}(\cite{gap-EJ} )\label{Projection}
	Given two  independent random variables $X, Y$ with absolutely continuous densities, and a function $f$ meeting
	$$\E[f^2(X+Y)]<\infty,\;\; \E[f(X+Y)]=0.$$
	Define
	$$g_1(u):=\E[f(u+Y)],\;\; g_2(v):=\E[  f(X+v)].$$
	Then for any measurable functions $h_1, h_2$,
	$$\E[f(X+Y)-h_1(X)-h_2(Y)]^2=\E[f(X+Y)-g_1(X)-g_2(Y)]^2+\E[g_1(X)-h_1(X)]^2+\E[g_2(Y)-h_2(Y)]^2.$$
\end{lemma}
The vertical component of the projection is control by Poincar$\acute{\mbox{e}}$ constant, provided that the expectation of random variables are all $0$ \cite{gap-EJ}, thus in Section \ref{proof of independent convolution J}, \textbf{we suppose all  random variables have 0 expectation without specifically statement.}

\begin{lemma}(\cite{gap-EJ} )\label{low bound of Poincare}
		Given two  independent random variables $X, Y$ with absolutely continuous densities, and a function $f$ meeting
	$$\E[f^2(X+Y)]<\infty,\;\; \E[f(X+Y)]=0.$$
	Define
	$$g_1(u):=\E[f(u+Y)],\;\; g_2(v):=\E[  f(X+v)].$$
	Then for any $\beta\in[0,1]$,
	$$E[f(X+Y)-h_1(X)-h_2(Y)]^2\ge \dfrac{1}{\bar{J}}\left(\dfrac{\beta}{R_X^*}\E[ g_1(X)-\mu X]^2+\dfrac{1-\beta}{R_Y^*}\E[g_2(Y)-\mu Y]^2  \right),$$
	where $\bar{J}=(1-\beta)J(X)+\beta J(Y), \mu=-2J(X+Y).$
	
\end{lemma}

The next lemma gives the calculation about  the score function of $X+Y$.

\begin{lemma}(\cite{gap-EJ} )\label{score of sum}
	Given two  independent random variables $X, Y$ with absolutely continuous densities,  then
	$$\rho_{X+Y}(x)=\E[\rho_X(X)|X+Y=x]=\E[\rho_Y(Y)|X+Y=x].$$
	Further,
	$$\dfrac{J(X)+J(Y)}{2}-J(\dfrac{X+Y}{\sqrt{2}})=2\E[\rho_{X+Y}(X+Y)-\dfrac{\rho_X(X)+\rho_Y(Y)}{2}  ]^2.$$
\end{lemma}

Finally, some calculations are given:
\begin{lemma}(\cite{gap-EJ} )\label{compute}
	Given two  independent random variables $X, Y$ with absolutely continuous densities, 
	define
	$$f=2\rho_{X+Y},\;\;  g_1(u):=\E[f(u+Y)].$$
	Then,
	$$\E[X\rho_X(X)]=\frac{1}{2}\E[g_1(X)X]+\frac{1}{2}\E[g_2(Y)Y]=-1.$$
	$$\E[\rho_{X+Y}(X+Y)\rho_X(X)]=J(X+Y).$$
	$$\E[X+\rho_X(X)]^2=J(X)-2+\E[X^2].$$
	$$\E[\rho_X(X)]=\E[g_1(X)]=0.$$
\end{lemma}

\subsection{Proof of \eqref{simple result like corollary}}
	
		\begin{theorem}\label{change of J}
		Given two  independent random variables $X, Y$ with absolutely continuous densities, and constant $R$ meeting $R\ge \max\{ R_X^*, R_Y^*\}$, then
		$$J(\dfrac{X+Y}{\sqrt{2}})- \dfrac{1}{\sigma^2}\le \dfrac{2R}{\sigma^2+2R} \dfrac{1}{2}[ J(X)+J(Y)-\dfrac{2}{\sigma^2}],$$
		where $\sigma^2=\frac{\sigma_X^2+\sigma_Y^2}{2}=\sigma^2_{\frac{X+Y}{\sqrt{2}}}.$
	\end{theorem}
	\pf
	
	We would give the estimations of $\E[\rho_{\frac{X+Y}{\sqrt{2}}}(\frac{X+Y}{\sqrt{2}})-\frac{g_1(x)+g_2(y)}{\sqrt{2}} ]^2$, $\E[g_1(X)+X]^2$ and $\E[g_2(Y)+Y]^2$.
	
	 On the one hand, due to Lemma \ref{Projection}, $\frac{g_1(X)+g_2(Y)}{\sqrt{2}}$ is the projection of $\rho_{\frac{X+Y}{\sqrt{2}}}(\frac{X+Y}{\sqrt{2}})$ on the space constructed by $X, Y$(see Fig.1 in \cite{gap-EJ}, page 400). On the other hand, denote $K$ as the pedal of $\rho_{\frac{X+Y}{\sqrt{2}}}(\frac{X+Y}{\sqrt{2}})$ on the line $\left( \frac{\rho_X(X)+\rho_Y(Y)}{\sqrt{2}}, -\frac{ X+Y}{\sqrt{2}} \right)$, thus
	\begin{equation}\label{K be low bound}
	\E[\frac{g_1(X)+g_2(Y)}{\sqrt{2}}+\frac{ X+Y}{\sqrt{2}} ]^2\ge \E[K+\frac{ X+Y}{\sqrt{2}}]^2,
	\end{equation}
	since the hypotenuse of a triangle is the longest.

	Note that $K$ can be expressed as $\lambda[\frac{\rho_X(X)+\rho_Y(Y)}{\sqrt{2}}]-(1-\lambda)[\frac{ X+Y}{\sqrt{2}}],$ thus
	$$\lambda=\arg\min_{a\in[0,1]} \E\left[ \rho_{\frac{X+Y}{\sqrt{2}}}(\frac{X+Y}{\sqrt{2}})-\left( a[\frac{\rho_X(X)+\rho_Y(Y)}{\sqrt{2}}]-(1-a)[\frac{ X+Y}{\sqrt{2}}] \right)    \right]^2.$$
	According to
	$$\arg\min_{a}\E[U-aV]^2=\dfrac{\E(UV)}{\E(V^2)},$$
	and Lemma \ref{compute}, 
	\begin{align}\label{A' and A}
	\lambda&=\dfrac{\E\left([ \rho_{\frac{X+Y}{\sqrt{2}}}(\frac{X+Y}{\sqrt{2}})+ \frac{ X+Y}{\sqrt{2}}][\frac{\rho_X(X)+\rho_Y(Y)}{\sqrt{2}}+\frac{ X+Y}{\sqrt{2}}]\right)}{\E\left(  \frac{\rho_X(X)+\rho_Y(Y)}{\sqrt{2}}+\frac{ X+Y}{\sqrt{2}}  \right)^2}\\
	&=\dfrac{J(\frac{ X+Y}{\sqrt{2}})-2+(\sigma_X^2+\sigma_Y^2)/2}{1/2J(X)+1/2J(Y)-2+(\sigma_X^2+\sigma_Y^2)/2}\\
	&\overset{\Delta}{=}\dfrac{A'}{A}.
	\end{align}
	Put \eqref{A' and A} in  \eqref{K be low bound},
	\begin{equation}\label{put in lambda}
	\E[\frac{g_1(X)+g_2(Y)}{\sqrt{2}}+\frac{ X+Y}{\sqrt{2}} ]^2=\dfrac{1}{2}\E[g_1(X)+X]^2 +\dfrac{1}{2}\E[g_2(Y)+Y]^2\ge \dfrac{A'^2}{A}.
	\end{equation}
	
	Next, define
	$$ h_1(x)=h_2(x)=-x.$$
	By Lemma \ref{Projection} and \eqref{put in lambda}
	$$\E \left( \rho_{\frac{X+Y}{\sqrt{2}}}(\frac{X+Y}{\sqrt{2}})-  \frac{g_1(X)+g_2(Y)}{\sqrt{2}}  \right)^2\le A'-\dfrac{A'^2}{A}.$$
	So due to Lemma \ref{low bound of Poincare}
	\begin{align*}
	&A'-\dfrac{A'^2}{A}\ge \E[\rho_{\frac{X+Y}{\sqrt{2}}}(\frac{X+Y}{\sqrt{2}})-  \frac{g_1(X)+g_2(Y)}{\sqrt{2}}  ]^2\\
	&\ge \dfrac{1}{2\bar{J}}\left(\dfrac{\beta}{R_X^*}\E[ g_1(X)-\mu X]^2+\dfrac{1-\beta}{R_Y^*}\E[g_2(Y)-\mu Y]^2  \right)\\
	&=\dfrac{1}{2\bar{J}}[\dfrac{\beta}{R_X^*}[\E(g_1(X)+X)^2-2(\mu+1)(\sigma_X^2-1)+\sigma_X^2(\mu+1)^2]\\
	&\;\;\;\;+\dfrac{1-\beta}{R_Y^*}[\E(g_2(Y)+Y)^2-2(\mu+1)(\sigma_Y^2-1)+\sigma_Y^2(\mu+1)^2] ].
	\end{align*}
	Take $\beta=\frac{R_X^*}{R_X^*+R_Y^*}$, and use \eqref{put in lambda} again,
	\begin{equation}\label{fundamantal inequality}
	\begin{aligned}
	&A'-\dfrac{A'^2}{A}\\
	&\ge \dfrac{1}{2\bar{J}}\left(\dfrac{\E[g_1(X)+X]^2+\E[g_2(Y)+Y]^2-2(\sigma_X^2+\sigma_Y^2-2)(1-J(\frac{X+Y}{\sqrt{2}})) + (\sigma_X^2+\sigma_Y^2)(1-J(\frac{X+Y}{\sqrt{2}}))^2}{R_X^*+R_Y^*}\right)\\
	&\ge \dfrac{1}{R_Y^*J(X)+R_X^*J(Y)}[\dfrac{A'^2}{A}-(\sigma_X^2+\sigma_Y^2-2)(1-J(\frac{X+Y}{\sqrt{2}}))+(J(\frac{X+Y}{\sqrt{2}})-1)^2\dfrac{\sigma_X^2+\sigma_Y^2}{2}   ].
	\end{aligned}
	\end{equation}
	
	Suppose $\dfrac{\sigma_X^2+\sigma_Y^2}{2}=1$, then
	$$A'-\dfrac{A'^2}{A}\ge \dfrac{1}{R_Y^*J(X)+R_X^*J(Y)}[\dfrac{A'^2}{A}+A'^2].$$
	Due to $R\ge \max\{ R_X^*, R_Y^*\}$, then
	$$A'-\dfrac{A'^2}{A}\ge \dfrac{1}{R[J(X)+J(Y)]}[\dfrac{A'^2}{A}+A'^2]=\dfrac{A'^2}{2R A}.$$
	Thus 
	$$A'\le \dfrac{2R}{1+2R}A,$$
	\textit{i.e.}
	$$J(\dfrac{X+Y}{\sqrt{2}})\le \dfrac{2R}{1+2R}[\dfrac{1}{2}(J(X)-1)+\dfrac{1}{2}(J(Y)-1)]+1.$$
	
	For general case $1/2\sigma_X^2+1/2\sigma_Y^2=\sigma^2>0$, 
	\begin{align*}
	J(\dfrac{X+Y}{\sqrt{2}})=\dfrac{1}{\sigma^2}J(\dfrac{X+Y}{\sigma\sqrt{2}}) \le \dfrac{1}{\sigma^2}\left( \dfrac{2R}{\sigma^2+2R}[ \dfrac{1}{2}(\sigma^2 J(X)-1)  +\dfrac{1}{2}(\sigma^2J(Y)-1)]+1\right),
	\end{align*}
\textit{	i.e. }
	$$J(\dfrac{X+Y}{\sqrt{2}})- \dfrac{1}{\sigma^2_{\frac{X+Y}{\sqrt{2}}}}\le \dfrac{2R}{\sigma^2+2R} \dfrac{1}{2}[ J(X)+J(Y)-\dfrac{4}{\sigma_X^2+\sigma_Y^2}].$$
	Then the first inequility has been proved. Recall $J_2\sigma^2\ge 1$, the second inequility is easy to obtain.
	
	\e

\end{document}